\documentclass[preprint,1p,12pt]{elsarticle}
\journal{Discrete Mathematics}
\usepackage[utf8]{inputenc}
\usepackage{mathrsfs}
\usepackage{enumitem}
\usepackage{complexity}
\usepackage{todonotes}
\usepackage{hyperref}
\usepackage{comment}
\usepackage[framemethod=TikZ]{mdframed}

\usepackage{amsfonts,amsmath,amssymb}

\tikzstyle{vertex}=[circle, draw, inner sep=0pt, minimum size=21pt]
\tikzstyle{svertex}=[circle, draw, inner sep=0pt, minimum size=3pt]
\tikzstyle{dvertex}=[circle, draw, inner sep=0pt, minimum size=9pt]
\tikzstyle{vertbox}=[draw, inner sep=3pt, minimum size=8pt]

\newtheorem{thm}{Theorem}
\newtheorem{lem}
{Lemma}
\newdefinition{rmk}{Remark}\newtheorem{proposition}{Proposition}
\newtheorem{cor}{Corollary}
\newproof{pf}{Proof}

\def\equationautorefname~#1\null{Equation~(#1)\null}
\def\thmautorefname~#1\null{Theorem~#1\null}
\def\lemautorefname~#1\null{Lemma~#1\null}\def\corautorefname~#1\null{Corollary~#1\null}\def\propositionautorefname~#1\null{Proposition~#1\null}\def\rmkautorefname~#1\null{Remark~#1\null}
\def\figureautorefname~#1\null{Figure~#1\null}
\def\sectionautorefname~#1\null{Section~#1\null}

\newcommand{\Oh}{\mathscr{O}}
\newcommand{\iso}{\mathsf{iso}}

\hypersetup{
	colorlinks,
	linkcolor=red!60!black,
	citecolor=blue!60!black,
	urlcolor=green!60!black
}

\begin{document}
\begin{frontmatter}
\title{Sum Labelling Graphs of Maximum Degree Two}
\author[1]{Henning Fernau (Universit\"at Trier, Universitätsring 15, Trier, Germany)%
}
\ead{fernau@uni-trier.de}

\author[2]{Kshitij Gajjar (Indian Institute of Technology Jodhpur, Rajasthan, India)
}
\ead{kshitij@iitj.ac.in}

%
%
%
    
%
\begin{abstract}

The concept of \emph{sum labelling} was introduced in 1990 by Harary. 
A graph is a~\emph{sum graph} if its vertices can be labelled by distinct positive integers in such a way that two vertices are connected by an edge if and only if the sum of their labels is the label of another vertex in the graph. It is easy to see that every sum graph has at least one isolated vertex, and every graph can be made a sum graph by adding at most $n^2$ isolated vertices to it. The minimum number of isolated vertices that need to be added to a graph to make it a sum graph is called the~\emph{sum number} of the graph.

The sum number of several prominent graph classes (e.g., cycles, trees, complete graphs) is already well known. We examine the effect of taking the disjoint union of graphs on the sum number. In particular, we provide a complete characterization of the sum number of graphs of maximum degree two, since every such graph is the disjoint union of paths and cycles.
\end{abstract}

\end{frontmatter}

\section{Introduction}

The area of graph labelling 
is a specific subarea of graph theory that has developed an enormous body of literature, as testified by Gallian's dynamic survey~\cite{Gal2020} which mentions over 3000 research papers.
One of these labellings is \emph{sum labelling}, introduced by Harary~\cite{harary} as a form of representing graphs. It is known~\cite{KraMilNgu2001} that every $n$-vertex graph~$G$ can be represented via a sum labelling, which means that it is possible to add at most $n^2$ isolated vertices (also called isolates, in short) to~$G$ to make it a sum graph. This makes sum labelling a compelling concept from the viewpoint of computer science also, because it may be that certain graphs can be encoded much more succinctly with sum labellings than with the more traditional ways of storing graphs.



\begin{figure}
\begin{centering}
\begin{tikzpicture}

\def \xoff {5};
\def \isoso {1};

\node[vertex](a1) at (-4, 4) {};
\node[vertex](a2) at (-3, 3) {};
\node[vertex](a3) at (-3, 5) {};

\node[vertex](b1) at (-4+\xoff,4) {$1$};
\node[vertex](b2) at (-3+\xoff,3) {$3$};
\node[vertex](b3) at (-3+\xoff,5) {$2$};

\node[vertex](iso1b) at (-3+\xoff+\isoso,3.4) {$5$};
\node[vertex](iso2b) at (-3+\xoff+\isoso,4.6) {$4$};

\node[vertex](d1) at (-4+2*\xoff,4) {$1$};
\node[vertex](d2) at (-3+2*\xoff,3) {$4$};
\node[vertex](d3) at (-3+2*\xoff,5) {$3$};

\node[vertex](iso1d) at (-3+2*\xoff+\isoso,3.4) {$7$};
\node[vertex](iso2d) at (-3+2*\xoff+\isoso,4.6) {$5$};

\node[] at (-3,2) {\Large (a)};
\node[] at (-3+\xoff,2) {\Large (b)};
\node[] at (-3+2*\xoff,2) {\Large (c)};

\begin{scope}[every path/.style={-}, every node/.style={inner sep=1pt}]
       \draw (a1) -- (a2) -- (a3) -- (a1);
       \draw (b1) -- (b2) -- (b3) -- (b1);
       \draw (d1) -- (d2) -- (d3) -- (d1);
       
\end{scope} 
\end{tikzpicture}
\caption{(a) This graph is not a sum graph, because it has no isolated vertices; (b) This is an incorrect sum labelling of a sum graph, because $(1,4)$ is not an edge yet there is a vertex labelled $1+4=5$ in the graph; (c) This is a correct sum labelling of a sum graph.}
\label{fig:validinvalidsumlabelling}
\end{centering}
\end{figure}
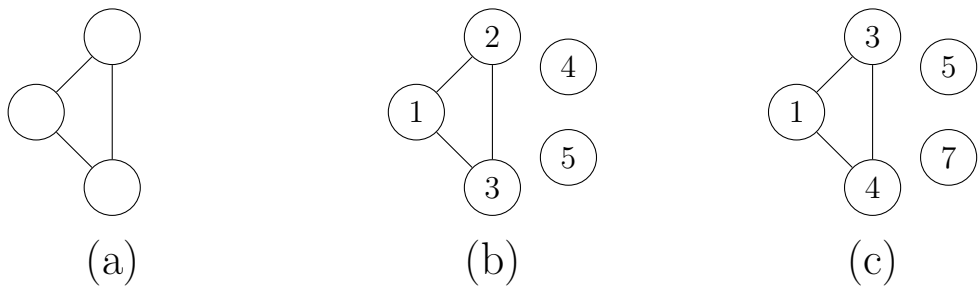

Let us now fix some notations. We deal with simple, undirected graphs, specified (as usual) as $G=(V,E)$, where $V$ is the (finite) set of vertices of $G$, and $E$ is its set of edges. If $v$ is an endpoint of an edge~$e$, then we say that $v$ and~$e$ are incident. The number of edges incident to a vertex is the degree of the vertex. Let $\mathbb{N}$ denote the set of all natural numbers (positive integers). Then, we say that $G$ is a \emph{sum graph} if there exists an injective mapping $\lambda:V\to\mathbb{N}$ (called the~\emph{sum labelling} of the vertices of~$G$) such that

$$E=\{xy\mid \exists z\in V:\lambda(z)=\lambda(x)+\lambda(y)\}.$$

Up to isomorphism, the set of numbers $\lambda(V)$ therefore determines~$G$. In other words, $\lambda$ encodes~$G$. As isolated vertices (i.e., vertices of degree zero) are usually irrelevant in applications, $\lambda(V)$ can be viewed as the description of $G\setminus I$, where $I$ is the set of all isolated vertices of $G$. Then, $\lambda(V)$ is called the~\emph{sum number encoding} of $G\setminus I$. Conversely, given a graph $G$ without isolates, the minimum number of isolates that need to be added to $G$ in order to make it a sum graph is called the \emph{sum number} of~$G$, written as $\sigma(G)$. Thus, $G+N_{\sigma(G)}$ is a sum graph. (Here, $+$ denotes the disjoint union of graphs. Also, $N_i$ denotes the null graph (edgeless graph) on $i$ vertices, or equivalently, a set of $i$ isolated vertices.) See~\autoref{fig:validinvalidsumlabelling} for some examples and non-examples of sum graphs and sum labellings.

A labelling function $\lambda$ can be also seen as operating on edges by the summability condition. $\lambda(e)$ for an edge $e=xy\in E$ is defined as $\lambda(x)+\lambda(y)$. Thus, though only the vertices are labelled by a sum-labelling, we sometimes also refer to its edges as labelled by the sum of its endpoints (two different edges can have the same edge label).

Are substantial savings possible with sum number encodings of graphs?
Some partial answers are possible from the literature. For instance, 
%
%
$\sigma(K_n)=2n-3$ is known for $n\geq 4$, i.e., $3n-3$ numbers suffice to store the information about the complete graph $K_n$, while 
traditional methods would need $\Oh(n^2)$ bits. As mentioned in~\cite{Smy91}, this can be obtained by labelling vertex $x_i$ with $4i-3$, with $1\leq i\leq n$, leading to isolate labels $4j+2$ for $1\leq j\leq  2n-3$. Hence, the sizes of the labels are in fact linear in~$n$. 

The focus of our study is the sum number of certain graphs. This follows much of the tradition in the literature, as can be seen in surveys like~\cite{Gal2020,Rya2009}. More precisely, we prove as our main result a complete picture of the sum number of every graph of maximum degree two. As a consequence, if $G$ has maximum degree two, then $\sigma(G)\leq 3$. This is not completely expected, as it is known that the sum number of  general graphs grows with the number of edges~\cite{NagMilSla2001}.
In fact, this can happen even with sparse graphs~\cite{HarSmy95,SutMil01}.


When talking about sum labelling a whole infinite family of graphs~$\cal G$, often with the additional property that for each positive integer~$n$, there is at most one graph $G_n$ of order $n$ within~$\cal G$, we also speak of a \emph{labelling scheme} $\lambda:\mathbb{N}\to\mathbb{N}$ that formalizes the labelling strategy that we suggest for~$\cal G$ in the following sense. For $G_n$, to be labelled with $i$ isolates, we take $\{1,\dots,n\}$ as the vertex set of $G_n$ and consider the set of numbers $\{\lambda(1),\dots,\lambda(n+i)\}$ as the set of labels of the sum graph $G_n+N_i$. 
Extending this notion, a \emph{general labelling scheme} is specified by three functions $\lambda:\mathbb{N}\to\mathbb{N}$,  $\sigma:\mathbb{N}\to\mathbb{N}$  and  $\iota:\mathbb{N}\to\mathbb{N}$  that are interpreted as a labelling strategy for $G_n\in\cal G$ of order~$n$, with $i$ isolates, as follows. As the set of vertices of  $G_n+N_i$, we consider $V_{n+i}=\{\sigma(n),\sigma(n)+1,\dots,\sigma(n)+n-1,\iota(n),\iota(n)+1,\dots, \iota(n)+i-1\}$, where the first $n$ numbers denote the vertices of~$G_n$, and as labels we take $\lambda(j)$ with $j\in V_{n+i}$. This boils down to a labelling scheme if $\sigma(n)$ is constant one and $\iota(n)=n+1$. More general labelling strategies of $n$-vertex graphs of a family of graphs~$\cal G$ are possible and will be discussed later in this paper.

Our main result is a complete precise characterization of all graphs $G$ of maximum degree two:

\begin{thm}\label{thm:main-maxdegtwo}
Let $G$ be a graph of maximum degree two. Then, $\sigma(G)=\delta(G)$ except for two graphs, namely $C_4$ and $C_4+P_2$, for which $\sigma(G)=\delta(G)+1$.
\end{thm}

Harary~\cite{harary} already showed that $\sigma(C_4)=3$, and that the minimum degree of a graph is always a lower bound on its sum number (i.e., $\sigma(G)\geq\delta(G)$ for all graphs $G$). Therefore, to prove our main theorem, it suffices to show that $\sigma(G)\leq\delta(G)$ for all graphs $G$ of maximum degree two, except for $C_4$ and $C_4+P_2$.
An additional proof is required to show that $\sigma(C_4+P_2)=2$.
 Apart from having a combinatorial result, we can also interpret our proof as providing an algorithm that labels any graph of maximum degree two optimally with respect to its sum number.

For the motivation of efficiently storing graphs, this is not completely satisfying, as the sizes of the labels could be exponential in the number of vertices of the graph according to our constructions, which means that we might need up to $\Oh(n^2)$ many bits for storing an $n$-vertex graph. In principle and in general, we can do this more efficiently in terms of label sizes~\cite{FerGaj2021}, but the algorithm presented in~\cite{FerGaj2021} is not tailored towards using as few isolates as possible, i.e., it does not obey the sum number of the graph, which is the focus of this study.


Notice that every graph of maximum degree two is a disjoint union of cycles and paths (in other words, each connected component of the graph is either a path or a cycle). To prove our main theorem, we will deal with the connected components in a specific sequence. This naturally produces an algorithm that optimally labels (with respect to the sum number) all graphs with maximum degree two. We provide a sketch of our strategy in~\autoref{fig:labelling-algo}.

\begin{figure}
\centering
\begin{mdframed}[backgroundcolor=blue!5, roundcorner=5pt]\begin{center}\underline{\bf Sum-labelling graphs of maximum degree two: a strategy}\end{center}
\begin{enumerate}
    \item Firstly, we deal with all cycles of length not equal to four (if any), in descending order of length.
    \item Secondly, we deal with all  cycles of length four (if any).
    \item Finally, we deal with paths (if any), in descending order of length.
\end{enumerate}
\end{mdframed}
\caption{Our proposed strategy for sum-labelling graphs $G$ with $1\leq\delta(G)\leq\Delta(G)\leq 2$.} \label{fig:labelling-algo}
\end{figure}

\autoref{fig:labelling-algo} also explains the sequence in which we will treat all graphs of maximum degree two. For example, if the graph $G$ is $$G = 5C_3 + 2C_4 + C_6 + 2C_7 + 3C_9 + 4P_2 + P_5 + 2P_8 + P_9,$$ then we will deal with the components of $G$ in the following order: $$3C_9,2C_7,C_6,5C_3,2C_4,P_9,2P_8,P_5,4P_2.$$

\section{The space complexity of sum labelling}
One of our motivations to return to sum labellings was the idea that one can use them to efficiently store graphs. This idea was already expressed in~\cite{KraMilNgu2001}. There, they consider the notion of the~\emph{range} $r(\lambda)$ of a labelling $\lambda$, which is defined as the difference between $\max \lambda(V)$ and $\min\lambda(V)$,\footnote{\label{fn:spum}
In~\cite{KraMilNgu2001} and also in~\cite{SinTiwTri2021}, under the name \emph{spum}, the mentioned difference is considered only for labellings that attain the sum number.} with $$r(\lambda)=\max_{v,v'\in V}\lambda(v)-\lambda(v')\,.$$ To clearly distinguish our notion of range from the ones mentioned in \autoref{fn:spum}, let us introduce the {\it sum range number} $r_\sigma(G)$ of a graph $G$ as the smallest range of a labelling of a sum graph $G+N_k$ \emph{for some~$k$}.
As eventually the range grows with the number of vertices, here we propose two different ways of ensuring that the numbers involved do not grow too fast. 

To better motivate the introduction of these new graph parameters, let us first analyze the sizes needed to store graphs in a database using a sum labelling encoding.
A graph~$G=(V,E)$ on $n$ vertices can be stored as follows: We need $\Oh(\log n)$ bits to store $n$ itself, plus  $\Oh(\log \log(\max\lambda(V)))$ bits to store $\log_2(\max\lambda(V))$, $\Oh(\log \sigma(G))$ bits to store the number of isolates and then $\log_2(2\max\lambda(V))\cdot (n+\sigma(G))$ more bits for the (at best ordered) list of numbers (vertex and isolate labels).
In the end, we have to store a list of $n+\sigma(G)$ many integers, each with $\log_2(2k)$ many bits, because edge labels (e.g., labels of isolates) have value of at most $2k$.

Instead, one could also first store the smallest label and then one would only need $\log_2(r)$ bits per number, where $r$ is the range of the labelling. More precisely, if we want to given an estimate of the number of bits needed to store graph $G=(V,E)$ with the labelling $\lambda$, we get the following formula.
\begin{equation}\label{eq:space-consumption}
    2(\log_2 n+\min_{v\in V}\log_2(\lambda(v)))+|\lambda(V\cup E)|\cdot \log_2(r(\lambda))
\end{equation}
Notice that although it looks beneficial to minimize $|\lambda(V\cup E)|$ by choosing a labelling $\lambda_\sigma$ that achieves $\sigma(G)$, i.e., where
$|\lambda_\sigma(V\cup E)|=|V|+\sigma(G)$, 
there could be another labelling $\lambda$ with  $|\lambda(V\cup E)|>|V|+\sigma(G)$, but $r(\lambda)$ could be much smaller than $r(\lambda_\sigma)$, potentially out-weighing the disadvantage of needing more isolates.
This is true in particular when $r$ takes values exponential in~$n$, as for the Ellingham-labelling for trees~\cite{ellingham}.

Further stretching our notation, we will also consider $r_\lambda$ for a labelling strategy $\lambda$, i.e., for a way to label a whole family of sparse graphs as described above, so that $r_\lambda$ can be viewed as a mapping that associates to $n$ the largest range of any labelling of an $n$-vertex graph according to this strategy. Hence, we can analyze the growth of $r_\lambda$ for certain labelling strategies.

What is the main purpose of a graph database?
Clearly, one has to access the graphs. A basic operation would be to answer the query if there is an edge between two vertices. 
Now, if $\max\lambda(V)$ is polynomial in $n=|V|$, we can answer this query 
in time $\Oh(\log(n))$. Namely, assuming the polynomial bound on the size of the labels, we would need time $\Oh(\log(n))$
to add the two labels of the vertices, and we also need time $\Oh(\log(n))$
to search for the sum in the ordered list of numbers, using binary search.
Otherwise, the additional time $\Oh(\log(\max\lambda(V)))$ would be quite expensive, probably making the idea of storing large graphs as sum graphs in databases unattractive. Therefore, also the range of labellings should be considered.

Other parameters that measure the space consumption of storing graphs even more accurately have been discussed in~\cite{FerGaj2021}.
However, for the discussions in this paper, the two parameters $\lambda$ and $r(\lambda)$ suffice, also because these are more accessible from the combinatorial viewpoint that we consider here.

The main difficulty in dealing with the combinatorics of sum labelling prevails also for these modified definitions, which is the question of how to prove lower bounds. The only general assertion that is available is to say that the sum number of a graph is at least as big as its minimum degree. There are also generalizations of this observation based on degree sequences (see \cite{Hao89,Smy91}), but this is irrelevant to us, as we consider graphs of bounded degree. For instance, this means that the sum number of a collection of cycles is at least two. But, as we see in the following, even proving that certain collections of cycles have a sum number of two is far from trivial. There are no really systematic tools available.

Regarding the notion of sum range number, it is nice to observe that the proof of Theorem~2.1 of \cite{SinTiwTri2021} concerning the spum of a graph is also valid in our case (which is, as discussed above, a definitorial variation of spum), so that we can state without proof the following result.

\begin{proposition}\label{prop:range-lb}
 Let $G$ be a graph of order $n$ with minimum degree $\delta(G)$ and maximum degree $\Delta(G)$. Then, $r_\sigma(G)\geq 2n-(\Delta(G)-\delta(G))-2$.   
\end{proposition}

Observe that for regular graphs, the lower bound stated in the previous proposition simplifies to $2n-2$. Unfortunately, even for our simple graph families, we reach this bound only occasionally.

\begin{figure}
\begin{centering}
\begin{tikzpicture}

\def \yoff {3.5};

\node[vertex](a2) at (-3.6,1.4) {$16$};
\node[vertex](a1) at (-3.1,0) {$31$};

\node[vertex](b1) at (-2.2,1.4) {$17$};
\node[vertex](b2) at (-1.7,0) {$30$};

\node[vertex](c1) at (-0.8,1.4) {$18$};
\node[vertex](c2) at (-0.3,0) {$29$};

\node[vertex](d1) at (0.6,1.4) {$19$};
\node[vertex](d2) at (1.1,0) {$28$};

\node[vertex](e1) at (2,1.4) {$20$};
\node[vertex](e2) at (2.5,0) {$27$};

\node[vertex](f1) at (3.4,1.4) {$21$};
\node[vertex](f2) at (3.9,0) {$26$};

\node[vertex](g1) at (4.8,1.4) {$22$};
\node[vertex](g2) at (5.3,0) {$25$};

\node[vertex](h1) at (6.2,1.4) {$23$};
\node[vertex](h2) at (6.7,0) {$24$};

\node[vertex](iso1) at (7.5,1.4) {$47$};

\node[] at (-4.8,0.7) {\large (b)};
\node[] at (-4.8,0.7+\yoff) {\large (a)};

\node[vertex](aa2) at (-3.6,1.4+\yoff) {$2$};
\node[vertex](aa1) at (-3.1,0+\yoff) {$3$};

\node[vertex](bb1) at (-2.2,1.4+\yoff) {$5$};
\node[vertex](bb2) at (-1.7,0+\yoff) {$6$};

\node[vertex](cc1) at (-0.8,1.4+\yoff) {$11$};
\node[vertex](cc2) at (-0.3,0+\yoff) {$12$};

\node[vertex](dd1) at (0.6,1.4+\yoff) {$23$};
\node[vertex](dd2) at (1.1,0+\yoff) {$24$};

\node[vertex](ee1) at (2,1.4+\yoff) {$47$};
\node[vertex](ee2) at (2.5,0+\yoff) {$48$};

\node[vertex](ff1) at (3.4,1.4+\yoff) {$95$};
\node[vertex](ff2) at (3.9,0+\yoff) {$96$};

\node[vertex](gg1) at (4.8,1.4+\yoff) {$191$};
\node[vertex](gg2) at (5.3,0+\yoff) {$192$};

\node[vertex](hh1) at (6.2,1.4+\yoff) {$383$};
\node[vertex](hh2) at (6.7,0+\yoff) {$384$};

\node[vertex](isoiso1) at (7.5,1.4+\yoff) {$767$};

\begin{scope}[every path/.style={-}, every node/.style={inner sep=1pt}]
       \draw (a1) -- (a2); \draw (aa1) -- (aa2);
       \draw (b1) -- (b2); \draw (bb1) -- (bb2);
       \draw (c1) -- (c2); \draw (cc1) -- (cc2);
       \draw (d1) -- (d2); \draw (dd1) -- (dd2);
       \draw (e1) -- (e2); \draw (ee1) -- (ee2);
       \draw (f1) -- (f2); \draw (ff1) -- (ff2);
       \draw (g1) -- (g2); \draw (gg1) -- (gg2);
       \draw (h1) -- (h2); \draw (hh1) -- (hh2);
       
\end{scope} 
\end{tikzpicture}
\caption{(a) The exponential labelling scheme; (b) The linear labelling scheme.}
\label{fig:disjointedges}
\end{centering}
\end{figure}
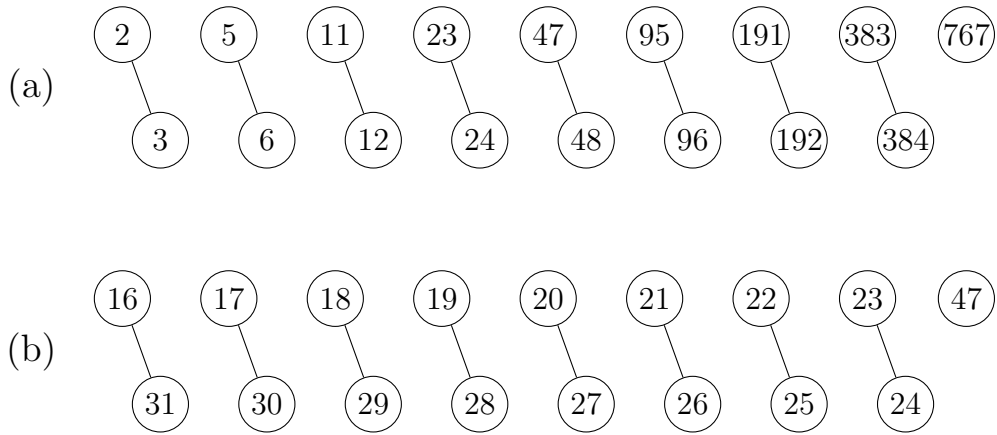

\section{A first example: labelling a disjoint collection of edges}


This section should be treated as an introductory example into the intricacies of sum labelling. It has also been studied earlier~\cite{FerGaj2021,SinTiwTri2021}.
Moreover, it covers an important subcase of our main theorem, which is 1-regular graphs, or graphs of (maximum) degree one (without isolates). Also, one can see examples that deal with the union of two graphs, each of sum number one.

It is known that all trees have sum number~1; according to a remark following Theorem 5.1 in~\cite{ellingham}, all forests also have sum number~1. However, it is not that clear how fast the label sizes grow in these constructions. 
Also, recall that it is still an open question for general graphs with sum number one whether their graph union again has sum number one~\cite{KraMilNgu2001}.
Thus, we will present two different constructions that label a disjoint collection of edges. More mathematically speaking, we will show two labelling schemes for the family of 1-regular graphs: an exponential labelling and a linear labelling. 

\subsection{An exponential solution}

If you have $n$ vertices (i.e., $n/2$ edges), label the first edge as
$(2,3)$. The second edge starts with the edge label of the first edge ($2+3=5$, so the second edge is labelled $(5,6)$). The third edge starts with the edge label of the second edge ($5+6=11$, so the third edge is labelled $(11,12)$), and so on (see~\autoref{fig:disjointedges} (a) for an example with $n=16$).

Generalising this, the following labelling scheme $\lambda:\mathbb{N}\to\mathbb{N}$ works for every 1-regular graph:
$$\lambda(n)=\left\{\begin{array}{ll}
     2&  \text{if }n=1\\
     \lambda(n-1)+1& \text{if $n$ is even}\\
     \lambda(n-2)+\lambda(n-1)& \text{if $n$ is odd and }n>1
\end{array}\right.$$

The Online Encyclopedia of Integer Sequences suggests that this is another variation on Ulam numbers if we think of the starting point to be $\lambda(0)=1$. Then, $\lambda(n)$ (for $n>1$) can be seen as the smallest (when $n$ is even) or largest (when  $n$ is odd)  number bigger than  $\lambda(n-1)$ that is a unique sum of two distinct earlier terms of the sequence. This connection also suggests the following closed form:
$$\lambda(n)=\left\{\begin{array}{ll}
     3\cdot 2^{k-1}& \text{if $n$ is even, i.e., }n=2k\\
     3\cdot 2^k-1& \text{if $n$ is odd, i.e., } n=2k+1
\end{array}\right.$$
In other words, we have $\lambda(n)\in\Theta\left(\left(\sqrt{2}\right)^n\right)$, implying that it is exponential in~$n$.
Although the suggested labelling~$\lambda$ is {optimal with respect to the sum number~$\sigma$}, we see: {$r(\lambda_n)\in\Theta\left(\left(\sqrt{2}\right)^n\right)$}. Can we do better with respect to the sum range number?

\subsection{A linear solution}

Consider the following general labelling scheme for 1-regular graphs (observe that~$n$ is necessarily even) that we first describe in a more intuitive fashion, already indicating the edges.
$$(n,2n-1),(n+1,2n-2),\ldots,\left(\frac{3n}{2}-1,\frac{3n}{2}\right)\,.$$ 
Here, writing $(\lambda(u),\lambda(v))$ refers to two vertices $u,v$ that are connected by an edge (see~\autoref{fig:disjointedges} (b) for an example with $n=16$). Notice that all edge labels sum to $3n-1$ (which is the isolate), and even the sum of the two smallest labels, i.e., $n+(n+1)=2n+1$, is smaller than $3n-1$ but bigger than any other label in the graph. More formally, we consider the functions $\lambda,\sigma,\iota$ with $\lambda(n)=\sigma(n)=n$ and $\iota(n)=3n-1$. This gives as vertex names $\{n,n+1,\dots,2n-1\}$ for a 1-regular graph of order~$n$.

This general labelling scheme can be further generalized by using the parameters $(x,y,d,k)$, with $x<y$ (in our example, $x=n,y=2n-1,d=1,k=n/2-1$), by putting $$(x,y),(x+d,y-d),\dots,(x+kd,y-kd)\,.$$
All labels sum up to $x+y$, which is the isolate.
As long as the sum of the two smallest labels, i.e., $2x+d$, is smaller than $x+y$ but bigger than $y$, such a sum labelling is valid. 
As the scheme consists of interleaving an increasing arithmetic progression with a decreasing arithmetic progression (with the same ``slope''), we call such schemes \emph{arithmetic progression schemes}.

The concrete arithmetic progression scheme that we first suggested has as its range the numbers $n$ through $2n-1$ and is hence (nearly) optimal, as \autoref{prop:range-lb} gives $2n-2$ as a lower bound. Singla, Tiwari \&  Tripathi~\cite{SinTiwTri2021} show an upper bound of $2n-1$. Therefore, we know that the optimal answer is either $2n-2$ or $2n-1$, but we do not know which one it is.





\subsection{Labelling paths}

The  ideas presented for 1-regular graphs work for paths also.
As we will need the exponential labelling scheme explicitly in the following, we are going to present (only) this one now.
For the linear solutions, we refer to \cite{FerGaj2021,SinTiwTri2021}.


A scheme could be based on fixing two positive integers $x,y$ as parameters, and then defining the labelling scheme $\lambda_{x,y}^\phi:\mathbb{N}\to\mathbb{N}$ as follows.

\begin{equation}\label{eq:fibo-scheme}
\lambda_{x,y}^\phi(n)=\left\{\begin{array}{ll}
     x&  \text{if }n=1\\
     y& \text{if } n=2\\
     \lambda_{x,y}^\phi(n-2)+\lambda_{x,y}^\phi(n-1)& \text{if  }n>2
\end{array}\right.
\end{equation}
Due to the similarity to Fibonacci numbers, it is clear that $\lambda_{x,y}^\phi(n)=\Oh(\phi^n)$, where $\phi$ is the golden ratio number, irrespectively of the start values $x,y$.
We can hence deduce the following well-known fact by this \emph{Fibonacci scheme}.
\begin{lem}For any $n\in\mathbb{N}$,
$\sigma(P_n)=1$.
\end{lem}

\section{Several labelling strategies for collections of cycles}


Recall that according to the algorithmic strategy sketched in~\autoref{fig:labelling-algo}, we first deal with all cycles of length five and larger, then with all triangles,  and finally with all cycles of length four.
The collection of $C_4$ is the most tricky one, as it could possibly leave us with three intermediate isolates.
Apart from this special situation, we will always face the situation that after having dealt with $k-1$ cycles, we have two isolates that we integrate into the $k$th cycle as the start of a new Fibonacci-type labelling. This is discussed in detail in the following subsections.

For the inductive argument, it becomes crucial to know that our labelling contains a \emph{non-trivial arithmetic progression}, or NTAP for short. This means that we find three labels $x,x+d,x+2d$ in the proposed labelling such that the \emph{offset}~$d$ is \emph{not} a label.

\subsection{Collections of 4-cycles}

In this subsection, we actually present two labelling strategies. The first one could be called ``linear-exponential'' in the sense that the proposed labelling strategy is linear (an arithmetic progression) per cycle, but from cycle to cycle, we observe an exponential growth. It uses three isolates (always) but has a smaller range compared to the second strategy that uses two isolates only (from two $C_4$ onwards) but needs a larger range.

\subsubsection{A linear-exponential labelling scheme}

Consider the labelling $(2,5,8,11)$ of a $C_4$.
Notice that the progression is arithmetic, with a difference of~3.
All numbers are congruent 2 modulo~3.

The three isolates are: $(7,13,19)$.
This arithmetic progression, with a difference of~6, can be again lifted to a labelling of a second $C_4$, which is then $(7,13,19,25)$. All numbers are congruent 1 modulo~3.

The three isolates are now: $(20,32,44)$. This arithmetic progression, with a difference of~12, can be again lifted to a labelling of a third $C_4$, which is then $(20,32,44,56)$. All numbers are congruent 2 modulo~3, as with the first $C_4$.

It is clear that we can continue this construction by adding a fourth $C_4$ with labels $(52,76,100,124)$. All numbers are congruent 1 modulo~3.

To wrap up, the odd-numbered cycles get numbers that are congruent to~1 modulo~3, while the even-numbered cycles get numbers which are congruent to~2 modulo~3.
These modulo~3 observations show that no edges can ever occur between vertices in subsequent cycles. As all the edge labels of the $i$th cycle can be found on the $(i+1)$th cycle, we can see that
(as the differences on the $i$th cycle are of the form $3\cdot 2^{i-1}$), the non-edges (diagonals) on the $i$th cycle cannot be represented by vertices on the $(i+1)$th cycle. 
By the aforementioned exponential growth of the labels one cycle to the next cycle, further non-edges cannot be represented by the suggested numbers. This proves:

\begin{lem}\label{lem:linear-exp:sum-bound}
If $G$ is a disjoint union of $C_4$'s, then $\sigma(G)\leq 3$.
\end{lem}

\noindent
Moreover, we can state:
\begin{lem}\label{lem:linear-exp:range-bound}
$r_\lambda(G)\in\Oh(2^{n/4})$ for a graph~$G$  of order~$n$ that is a union of~$C_4$'s, 
for the specific labelling scheme~$\lambda$ that we described above.
\end{lem}

\subsubsection{Towards optimal sum labellings}

We know that $\sigma(kC_4)\in\{2,3\}$ (\autoref{lem:linear-exp:sum-bound}), and it is known that $\sigma(C_4)=3$.
Can we possibly also show that $\sigma(2C_4)=2$ or even $\sigma(3C_4)=2$?
Let us try a bit of algebra, assuming arithmetic progression labellings of the two considered $C_4$'s.

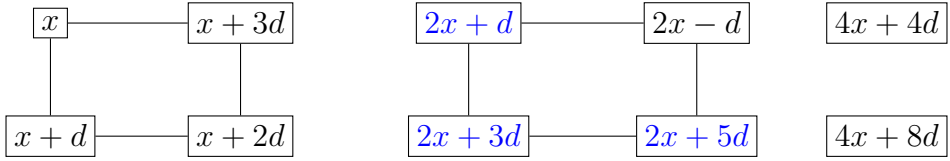
\begin{figure}[tbh]
\begin{centering}\scalebox{1}{\begin{tikzpicture}
\node[vertbox](a2) at (-4,1.5) {$x$};
\node[vertbox](a1) at (-4,0) {$x+d$};

\node[vertbox](b2) at (-1.5,1.5) {$x+3d$};
\node[vertbox](b1) at (-1.5,0) {$x+2d$};

\node[vertbox](d1) at (2-.5,1.5) {\color{blue}$2x+d$};
\node[vertbox](d2) at (2-.5,0) {\color{blue}$2x+3d$};

\node[vertbox](e1) at (4.5,1.5) {$2x-d$};
\node[vertbox](e2) at (4.5,0) {\color{blue}$2x+5d$};

\node[vertbox](iso1) at (7,1.5) {$4x+4d$};\node[vertbox](iso2) at (7,0) {$4x+8d$};

\begin{scope}[every path/.style={-}, every node/.style={inner sep=1pt}]
       \draw (a1) -- (a2); 
       \draw (b1) -- (b2); 
       \draw (a1) -- (b1); 
       \draw (d1) -- (d2); 
       \draw (e1) -- (e2); 
       \draw (a2) -- (b2); 
       \draw (d1) -- (e1); 
       \draw (d2) -- (e2); 
  
\end{scope} 
\end{tikzpicture}}
\caption{An algebraic approach to the $C_4$ problem.\label{fig:C4-algebra}}
\end{centering}
\end{figure}

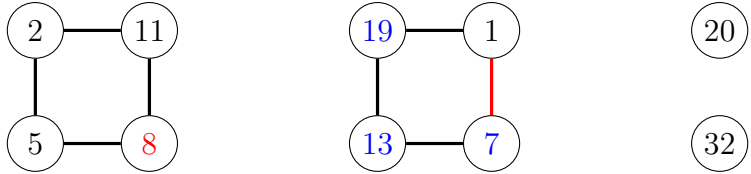
\begin{figure}[tbh]
\begin{centering}\scalebox{1}{

\begin{tikzpicture}

\node[vertex](a2) at (-3.0,1.5) {$2$};
\node[vertex](a1) at (-3,0) {$5$};
\node[vertex](b2) at (-1.5,1.5) {$11$};
\node[vertex](b1) at (-1.5,0) {\color{red}$8$};

\node[vertex](d1) at (1.5,1.5) {\color{blue}$19$};
\node[vertex](d2) at (1.5,0) {\color{blue}$13$};
\node[vertex](e1) at (3,1.5) {$1$};
\node[vertex](e2) at (3,0) {\color{blue}$7$};

\node[vertex](iso1) at (6,1.5) {$20$};\node[vertex](iso2) at (6,0) {$32$};


\begin{scope}[every path/.style={-}, every node/.style={inner sep=1pt}]
       \draw[very thick] (a1) -- (a2); 
       \draw[very thick] (b1) -- (b2); 
       \draw[very thick] (a1) -- (b1); 
       \draw[very thick] (d1) -- (d2); 
       \draw[very thick,red] (e1) -- (e2); 
       \draw[very thick] (a2) -- (b2); 
       \draw[very thick] (d1) -- (e1); 
       \draw[very thick] (d2) -- (e2); 
  
\end{scope} 
\end{tikzpicture}}
\caption{A minimal way to label a $2C_4$ with two isolates.\label{fig:2C4-labelling}}
\end{centering}
\end{figure}

The idea of~\autoref{fig:C4-algebra} is to find one of the three isolates of the second cycle within the labels of the first cycle. The only way this could happen is for the isolate $4x=(2x+d)+(2x-d)$.
Clearly, $4x\neq x$. If $4x=x+d$, then $3x=d$. This contradicts the label $2x-d$, which implies that $2x>d$.
If $4x=x+2d$, we conclude $3x=2d$, so that $x$ is even and $d$ is divisible by~$3$. The smallest numbers satisfying these conditions are $x=2$ and $d=3$; see \autoref{fig:2C4-labelling}. These divisibility conditions also enforce that all other labellings of this form have to be scalings of this minimal labelling by some constant factor. Finally, if $4x=x+3d$, then $x=d$. Hence, the number $2x+d=x+2d$ would occur twice as a vertex label. Therefore, under the conditions that our first cycle is labeled as in \autoref{fig:C4-algebra}, 
\autoref{fig:2C4-labelling}
basically shows the only possibility. Notice that this labelling contains the NTAP $2-5-8$.

\begin{figure}[tbh]
\begin{centering}\scalebox{.9}
{
\begin{tikzpicture}

\node[vertex](a2) at (-3.0,1.5) {\color{blue}$8$};
\node[vertex](a1) at (-3,0) {\color{green}$20$};
\node[vertex](b2) at (-1.5,1.5) {$44$};
\node[vertex](b1) at (-1.5,0) {\color{red}$32$};

\node[vertex](d1) at (0.5,1.5) {$28$};
\node[vertex](d2) at (0.5,0) {$52$};
\node[vertex](e1) at (2,1.5) {$4$};
\node[vertex](e2) at (2,0) {$76$};

\node[vertex](f1) at (4,1.5) {$5$};
\node[vertex](f2) at (4,0) {$3$};
\node[vertex](g1) at (5.5,1.5) {$15$};
\node[vertex](g2) at (5.5,0) {$17$};

\node[vertex](iso2) at (7.5,1.5) {$80$};\node[vertex](iso3) at (7.5,0) {$128$};

\begin{scope}[every path/.style={-}, every node/.style={inner sep=1pt}]
       \draw[very thick] (a1) -- (a2); 
       \draw[very thick] (b1) -- (b2); 
       \draw[very thick] (a1) -- (b1); 
       \draw[very thick] (d1) -- (d2); 
       \draw[very thick] (e1) -- (e2); 
       \draw[very thick] (a2) -- (b2); 
       \draw[very thick, red] (d1) -- (e1); 
       \draw[very thick] (d2) -- (e2); \draw[very thick, green] (g1) -- (f1); 
       \draw[very thick, green] (g2) -- (f2); 
       \draw[very thick, blue] (f1) -- (f2); 
       \draw[very thick, red] (g1) -- (g2); 
\end{scope} 
\end{tikzpicture}}\caption{\label{fig:3C4-labelling}A way to label a $3C_4$ with two isolates.}
\end{centering}
\end{figure}
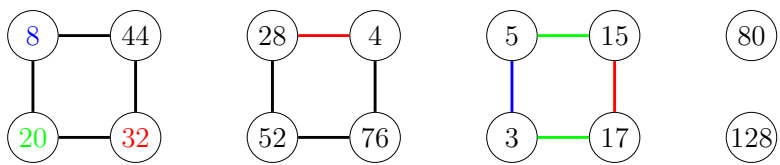

Could scaling help to also label $3C_4$ with our strategy? The somewhat surprising answer is \emph{yes}.
First, we look at a concrete example in \autoref{fig:3C4-labelling}.
The trick consists in the following steps:
\begin{enumerate}
    \item Multiply all labels used so far by a sufficiently large constant~$z>2$, which is four in our example. We actually need that (modulo~$z$) $z-1\neq z+1$. To ease our inductive argument, let us always pick $z=4$.
    \item Pick the smallest three labels of the first cycle, which is $x=8,x+d=20,x+2d=32$ in our example, and select numbers $a,b,c,e$ to label the third cycle. To avoid unwanted edges, choose $a=x/2+1$ (recall that $x$ must be an even number), $b=x/2-1$, $c=(x+2d)/2+1$, $e=(x+2d)/2-1$.
    \item Observe that the isolates of the $2C_4$-construction remain untouched.
    \item Also, since our labelling of $2C_4$ contains a NTAP, the proposed labelling of $3C_4$ contains a NTAP too.
\end{enumerate}

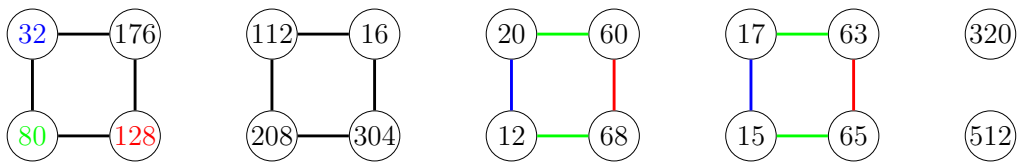
\begin{figure}[tbh]
\begin{centering}\scalebox{.9}{\begin{tikzpicture}
\node[vertex](a2) at (-3.0,1.5) {\color{blue}$32$};
\node[vertex](a1) at (-3,0) {\color{green}$80$};

\node[vertex](b2) at (-1.5,1.5) {$176$};
\node[vertex](b1) at (-1.5,0) {\color{red}$128$};

\node[vertex](d1) at (0.5,1.5) {$112$};
\node[vertex](d2) at (0.5,0) {$208$};

\node[vertex](e1) at (2,1.5) {$16$};
\node[vertex](e2) at (2,0) {$304$};

\node[vertex](f1) at (4,1.5) {$20$};
\node[vertex](f2) at (4,0) {$12$};
\node[vertex](g1) at (5.5,1.5) {$60$};
\node[vertex](g2) at (5.5,0) {$68$};

\node[vertex](h1) at (7.5,1.5) {$17$};
\node[vertex](h2) at (7.5,0) {$15$};

\node[vertex](i1) at (9,1.5) {$63$};
\node[vertex](i2) at (9,0) {$65$};

\node[vertex](iso2) at (11,1.5) {$320$};\node[vertex](iso3) at (11,0) {$512$};

\begin{scope}[every path/.style={-}, every node/.style={inner sep=1pt}]
       \draw[very thick] (a1) -- (a2); 
       \draw[very thick] (b1) -- (b2); 
       \draw[very thick] (a1) -- (b1); 
       \draw[very thick] (d1) -- (d2); 
       \draw[very thick] (e1) -- (e2); 
       \draw[very thick] (a2) -- (b2); 
       \draw[very thick] (d1) -- (e1); 
       \draw[very thick] (d2) -- (e2); \draw[very thick,green] (g1) -- (f1); 
       \draw[very thick,green] (g2) -- (f2); 
       \draw[very thick,blue] (f1) -- (f2); 
       \draw[very thick,red] (g1) -- (g2); 
       \draw[very thick,green] (h1) -- (i1); 
       \draw[very thick,green] (h2) -- (i2); 
       \draw[very thick,blue] (h1) -- (h2); 
       \draw[very thick,red] (i1) -- (i2); 
\end{scope} 
\end{tikzpicture}}\caption{\label{fig:4C4-labelling}A way to label a $4C_4$ with two isolates.}
\end{centering}
\end{figure}

As the $2C_4$-construction remains untouched up to scaling, we can actually repeat this argument, which could give the labelling of a $4C_4$ as in \autoref{fig:4C4-labelling}, and this type of argument continues to prove by induction on~$k$:

\begin{lem}\label{lem:kC4}
$\sigma(kC_4)=2$ for all $k\geq 2$. Moreover, the corresponding labelling contains a NTAP.
\end{lem}

\begin{pf}
Let us describe some details of the induction. For our inductive argument to work, we make the additional claim that the three smallest numbers $x,y,z$ labelling the first cycle form an arithmetic progression, i.e., there is a number $d$ such that $y=x+d$ and $z=x+2d$. Moreover, $d$ is not a label of any vertex, so that the labelling satisfies NTAP. 
The induction basis for $k=2$ was given above and satisfies NTAP. Let us assume that the labelling strategy works for some specific $k=K\geq 2$. When we multiply all labels of the first $K$ cycles by four, then this will not change the fact that (exactly) the edges of the $K$ cycles are described by these numbers, plus the two isolates that remain as isolates in the overall labelling. Also some NTAP is found after the modification by multiplication. The labelling of the $(K+1)$st cycle builds upon the smallest three labels $x,x+d,x+2d$ of the first cycle, choosing $a=x/2+1$, $b=x/2-1$, as well as $c=a+d$ and $e=b+d$ as labels of the last cycle. As we multiplied all original numbers by four, $x$ is an even number. Also, $a+b=x$, $a+e=b+c=x+d$ and $c+e=x+2d$, so that all wanted edge labels can be found as vertex labels on the first cycle. By way of contrast, the unwanted edges corresponding to $a+c=2a+d=x+d+2$ and $b+e=2b+d=x+d-2$ cannot be found as vertex labels, because all vertex labels of the first $K$ cycles (and also the isolates) are divisible by four, including the label $x+d$. Finally, as all `new labels' are  odd and all `old labels' are even, an edge between an `old vertex' and a `new vertex' must be labelled with a `new label', and this also implies that only the two bigger `new labels' $c$ and $e$ could possibly serve as edge labels. Moreover, as $c$ is one congruent four, this must match the only other label that is one congruent four, which is~$a$, as all `old labels' are divisible by four. Hence, the question is if $c-a=d$ is an `old label', which is clearly not the case by induction. Similarly, $b$ and $e$ are three congruent four, but $e-b=d$ and the same argument applies in this case as well.\qed 
\end{pf}

Notice that in the recursive labelling algorithm hidden in the previous proof, the assumption that $d$ does not occur as a vertex label is crucial, as otherwise there would be an unwanted edge between $a$ and $d$, because we have the vertex label $a+d$.

In contrast to the labelling strategy described in the previous subsection, and in particular analyzed in \autoref{lem:linear-exp:range-bound},
we obtain a worse relation concerning the growth of the range for this new labelling strategy.

\begin{lem} \label{lem:cycles-sum-opt:range-bound}
There is a labelling strategy $\lambda$ for disjoint unions of $C_4$'s such that $r_\lambda(G)\in \mathcal{O}(2^{n/2})$ for a graph of order~$n$ which is a collection of $n/4$ many $C_4$'s.
\end{lem}

\begin{pf}
Although also the size of the smallest label grows in this magnitude, it suffices to estimate the size of the largest label. For $2C_4$, this is $32=2\cdot 2^{8/2}$, as described in \autoref{fig:2C4-labelling}. As this largest label is always multiplied by $4=2^{4/2}$, we get $2\cdot 2^{n/2}=2\cdot 2^{(n-4)/2}\cdot 2^{4/2}$ by induction, considering an $n$-vertex graph which is a collection of $n/4$ many $C_4$'s.
\qed\end{pf}

\noindent
We can generalize the construction of \autoref{lem:kC4} to get the following result.

\begin{proposition}\label{prop:induction-step+C4s} Let $G$ be a graph with $\delta(G)=\sigma(G)=2$ such that there is a sum-optimal labelling $\lambda$ of the sum graph $H=G+N_2$ such that $\lambda(V(H))$ contains a NTAP,
then there is a sum-optimal labelling  $\lambda'$ of the sum graph $H'=G+C_4+N_2$ such that $\lambda(V(H'))$ contains a NTAP.
\end{proposition}


For example, consider $C_3+C_4$, starting with a labelling $1-3-4$ of the $C_3$. Now, the isolates are $5$ and $7$. Observe that $1,3,5$ is an arithmetic progression whose offset~$2$ is not a vertex label. Hence, we can multiply the numbers by~4 to get a labelling $4-12-16$ of the $C_3$, with isolate labels $20$ and~$28$. Finally, we label the $C_4$ as $1-3-9-11$. Clearly, the same isolate labels of $20$ and $28$ suffice. 

This proposition will come in handy when finally combining the results of this subsection with that of the next one.
The importance of the arithmetic progression becomes also clear when revisiting \autoref{lem-4cycle3ariprog}.
Unfortunately, as already described in \autoref{lem:cycles-sum-opt:range-bound}, the range will grow exponentially with base $\sqrt{2}$ if this proposition is applied repeatedly.

\subsection{Collections of cycles without 4-cycles}


We first discuss labelling strategies for collections of cycles without $C_4$, but \autoref{prop:induction-step+C4s} immediately shows a way how to add $C_4$ afterwards, as we explain in the following.

\subsubsection{Dealing with long cycles}
We will consider Fibonacci-labellings of cycles $C_n$, with $n> 4$,\footnote{Collections of such longer cycles were treated in~\cite{BurRusSug2008} in a similar fashion.}  leaving the case of triangles to be treated later. 
$$\lambda_n(x,y):(x,x+y,2x+y,3x+2y,\dots,F(n)x+F(n-1)y)$$
with isolates $ (F(n)+1)x+F(n-1)y$ and $F(n+1)x+F(n)y$,
using the 
Fibonacci sequence $F:(1,1,2,3,5,8,13,21,\dots)$ with $F(0)=0$ for convenience.
\begin{lem}
For any $n\geq 3$, $n> 4$, and any 
$1\leq x\leq y$,
$\lambda_n(x,y)$ gives a sum labelling of $C_n$.
\end{lem}
In fact, we can consider the labelling scheme $\lambda(x,y)(k)=F(k)x+F(k-1)y$.
Notice that the corresponding range function grows as $\varphi^n$ for $C_n$, where $\varphi$ is the golden ratio number.

The parameters $x$ and $y$ ($x<y$) give us great flexibility. We could start labelling the first of a certain number of longer cycles, starting with $x_1=1$ and $y_1=1$.
If the first cycle has length $n_1$, the last of its vertices would be labelled $F(n_1)x_1+F(n_1-1)y=F(n_1)+F(n_1-1)=F(n_1+1)$. The two isolates $ (F(n_1)+1)+F(n_1-1)=F(n_1+1)+1,F(n_1+1)+F(n_1)=F(n_1+2)$, giving us a new pair $(x_2,y_2)$. If the second cycle has length $n_2$,  the last of its vertices would be labelled $F(n_2)x_2+F(n_2-1)y_2=F(n_2)(F(n_1+1)+1)+F(n_2-1)F(n_1+2)$. We get similar expressions for the isolates, which will again form the start $(x_3,y_3)$ of labelling the next cycle etc.

We explain this labeling strategy by an example, a collection of three $C_5$ in \autoref{fig:three-C5}.
The disadvantage of this strategy comes from the fact that we do not see an arithmetic progression in it such that its offset is not a vertex label. In particular, how do we label $C_5+C_4$?

However, there is a remedy to it:
Singla, Tiwari \&  Tripathi~\cite{SinTiwTri2021} showed that (for $r_\sigma$, the spum number)
$r_\sigma(C_n)\in [2n-2,2n-1]$ for $n\geq 4$, and $r_\sigma(C_n)=2n-1$ for $n\geq 13$.
Namely, for odd $n\geq 5$, they propose the label set
$$L(C_n)=[n-3,2n-4]\cup \{3n-6,3n-4\}\qquad\text{ordered as}$$
$$(n-3,
{n-1},
{2n-5},
{n+1},
{2n-7},
{n+3},
{2n-9},\dots,2n-4,n-2,n-3)$$
For instance, this gives the labelling $(2,4,5,6,3)$ of a $C_5$, with isolates $9,11$. 
Unfortunately, this is not a valid labelling as claimed in the paper, as we get the unwanted edge between the vertices labelled~$2$ and~$9$, adding up to the label~$11$.
The labelling works for $n=7$, though, where we get $(4,6,9,8,7,10,5)$ with isolates $15,17$. Namely, we find that the difference $3n-6-(3n-4)=2$ between the two isolate labels only occurs in $[n-3,2n-4]$ when $n=5$.  
In~\cite{SinTiwTri2021}, another labelling was proposed for even $n\geq 4$:\footnote{Other than claimed, the proposed labelling actually does not work for $n=4$.} 
$$L(C_n)=[n-2,2n-3]\cup \{3n-5,3n-3\}\qquad\text{ordered like}$$
$$(n-2,{2n-3},n,2n-5,n+2,{2n-7},\dots,2n-4,n-1,n-2)$$
Thus, for $n=6$, the $C_6$ can be labelled $(4,9,6,7,8,5)$, with isolates $13,15$. In all cases, we clearly find a NTAP, for instance for the proposed $C_6$-labelling $4,5,6$ with offset~1. 

As this is of importance for our algorithm, let us show that there does indeed exist a sum labelling of the $C_5$ that satisfies the conditions we need.

\begin{lem}\label{lem:C5-special}
The labelling $(1,2,7,9,3)$, with isolates $4,12,16$ contains the arithmetic progression $2-7-12$ whose offset~$5$ is not a label of any vertex. 
\end{lem}

This looks worse than what the Fibonacci labelling would deliver, as we need three isolates, but it is in fact better, as $4,12,16$ might be seen as the start of a Fibonacci labelling on the second cycle. This proves that $\sigma(C_5+C_n)=2$ if $n\neq 4$, and accordingly an optimal labelling can be given that contains a NTAP. 
For instance, $C_5+C_3$ can be labelled with $\{1,2,3,4,7,9,12,16\}$, with isolates $\{20,28\}$.
%
%
%
Only for the next cycles, we employ the Fibonacci scheme. This preserves the property to have a NTAP.

\begin{lem}\label{lem:C5-C4-special}
The labelling of $C_5 + C_4$ with the $C_5$ labelled as $(2,4,6,10,16)$ and the $C_4$ labelled as $(1,9,18,26)$ (with isolates $27,44$) contains the NTAP
$10-27-44$. 
\end{lem}

\begin{pf}
The $C_5$-labelling follows a Fibonacci scheme. As $1+9=10$, one $C_4$ edge label is in the $C_5$, while its other edge labels are in the isolates. \qed
\end{pf}

This labelling scheme can be generalized as follows:
\begin{align*}
    C_5:&\ (a,b,a+b,a+2b,2a+3b);\\
    C_4:&\ (c,2b+c,3a+3b,3a+5b);\\
    \text{isolates}:&\ 3a+5b+c, 6a+8b.
\end{align*}

The above labelling requires that $a=2c$ (to ensure that $a+2b=2b+2c$ for the edge connecting $2b+c$ and $c$), and $b\neq 3c$ (to avoid the unwanted edge $a+(a+b)=b+4c$ being represented by $2b+c$).

\begin{figure}[bt]
\centering
	\begin{tikzpicture}[transform shape]
 \def\xoff{4}
			\tikzset{every node/.style={circle,minimum size=0.3cm}}
			\node[draw] (a) at (0,1) {{{1}}};
			\node[draw] (b) at (-0.95106,0.30902) {{{2}}};
			\node[draw] (c) at (-0.58779,-0.80902) {{{7}}};
			\node[draw] (d) at (0.58779,-0.80902) {{{9}}};
			\node[draw] (e) at (0.95106,0.30902) {{{3}}};
			\path (a) edge[-,very thick] (b);
			\path (b) edge[-,very thick] (c);
			\path (c) edge[green,-,very thick] (d);
			\path (d) edge[blue,-,very thick] (e);
			\path (e) edge[red,-,very thick] (a);
   \node[draw] (a1) at (0+\xoff,1) {{\textcolor{blue}{12}}};
			\node[draw] (b1) at (-0.95106+\xoff,0.30902) {{\textcolor{green}{16}}};
			\node[draw] (c1) at (-0.58779+\xoff,-0.80902) {{{28}}};
			\node[draw] (d1) at (0.58779+\xoff,-0.80902) {{{44}}};
			\node[draw] (e1) at (0.95106+\xoff,0.30902) {{\textcolor{red}{4}}};
			\path (a1) edge[-,very thick] (b1);
			\path (b1) edge[-,very thick] (c1);
			\path (c1) edge[purple,-,very thick] (d1);
			\path (d1) edge[orange,-,very thick] (e1);
			\path (e1) edge[-,very thick] (a1);
       
   \node[draw] (a2) at (0+2*\xoff,1) {{\textcolor{orange}{48}}};
			\node[draw] (b2) at (-0.95106+2*\xoff,0.30902) {{\textcolor{purple}{72}}};
			\node[draw] (c2) at (-0.58779+2*\xoff,-0.80902) {{{120}}};
			\node[draw] (d2) at (0.58779+2*\xoff,-0.80902) {{{192}}};
			\node[draw] (e2) at (0.95106+2*\xoff,0.30902) {{{312}}};
			\path (a2) edge[-,very thick] (b2);
			\path (b2) edge[-,very thick] (c2);
			\path (c2) edge[-,very thick] (d2);
			\path (d2) edge[-,very thick] (e2);
			\path (e2) edge[-,very thick] (a2);
        \end{tikzpicture}
        
        \caption{\label{fig:three-C5-with-Special}How to label a collection of cycles, here three $C_5$'s, with isolates $360$, $504$, and NTAP $2-7-12$.}
\end{figure}
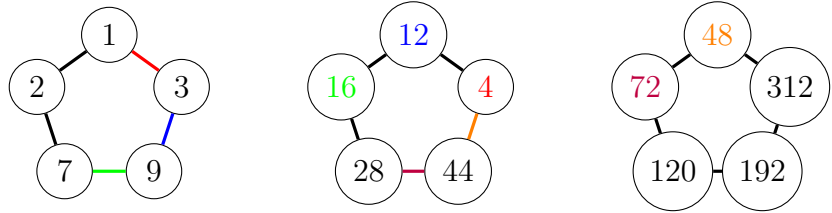

To conclude this section, we show how to label $3C_5$ in two different ways.
\autoref{fig:three-C5} shows a labelling were we strictly follow a Fibonacci scheme. 
If we use \autoref{lem:C5-special} at the beginning (with the advantage of showing an arithmetic progression as desired), we arrive at \autoref{fig:three-C5-with-Special}. 

\begin{figure}[bt]
\centering
	\begin{tikzpicture}[transform shape]
 \def\xoff{4}
			\tikzset{every node/.style={circle,minimum size=0.3cm}}
			\node[draw] (a) at (0,1) {{{1}}};
			\node[draw] (b) at (-0.95106,0.30902) {{{2}}};
			\node[draw] (c) at (-0.58779,-0.80902) {{{3}}};
			\node[draw] (d) at (0.58779,-0.80902) {{{5}}};
			\node[draw] (e) at (0.95106,0.30902) {{{8}}};
			\path (a) edge[-,very thick] (b);
			\path (b) edge[-,very thick] (c);
			\path (c) edge[-,very thick] (d);
			\path (d) edge[red,-,very thick] (e);
			\path (e) edge[blue,-,very thick] (a);
   \node[draw] (a1) at (0+\xoff,1) {{\textcolor{blue}{9}}};
			\node[draw] (b1) at (-0.95106+\xoff,0.30902) {{\textcolor{red}{13}}};
			\node[draw] (c1) at (-0.58779+\xoff,-0.80902) {{{22}}};
			\node[draw] (d1) at (0.58779+\xoff,-0.80902) {{{35}}};
			\node[draw] (e1) at (0.95106+\xoff,0.30902) {{{57}}};
			\path (a1) edge[-,very thick] (b1);
			\path (b1) edge[-,very thick] (c1);
			\path (c1) edge[-,very thick] (d1);
			\path (d1) edge[orange,-,very thick] (e1);
			\path (e1) edge[green,-,very thick] (a1);
       
   \node[draw] (a2) at (0+2*\xoff,1) {{\textcolor{green}{66}}};
			\node[draw] (b2) at (-0.95106+2*\xoff,0.30902) {{\textcolor{orange}{92}}};
			\node[draw] (c2) at (-0.58779+2*\xoff,-0.80902) {{{158}}};
			\node[draw] (d2) at (0.58779+2*\xoff,-0.80902) {{{250}}};
			\node[draw] (e2) at (0.95106+2*\xoff,0.30902) {{{408}}};
			\path (a2) edge[-,very thick] (b2);
			\path (b2) edge[-,very thick] (c2);
			\path (c2) edge[-,very thick] (d2);
			\path (d2) edge[-,very thick] (e2);
			\path (e2) edge[-,very thick] (a2);
        \end{tikzpicture}
        
        \caption{\label{fig:three-C5} How to label a collection of cycles, here three $C_5$'s, with isolates $474$, $658$.}
\end{figure}
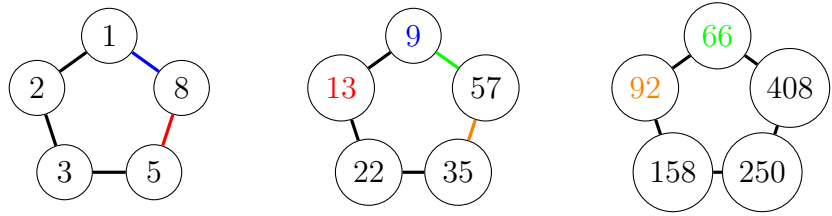

\subsubsection{Dealing with triangles}

We will prove the following assertion about collection of triangles ($C_3$'s).

\begin{lem}\label{lem:C3s}
Any Fibonacci labelling scheme for a non-empty collection of triangles gives a valid sum labelling that contains a NTAP.
\end{lem}


More precisely, consider $$\lambda_n'(x,y):(x,x+y,2x+y;3x+y,3x+2y,6x+3y;
9x+4y,9x+5y,18x+9y;\dots)$$
as a labelling scheme serving for $\ell$ many $C_3$'s, with isolates $$3^\ell x+\lfloor 3^\ell/2\rfloor y,\quad 3^\ell x+\lceil 3^\ell/2\rceil y\,.$$
This scheme is different in terms of growth compared to the schemes set up before for longer cycles.
Yet, it can be embedded in an inductive construction of a labelling of any number of cycles excluding $C_4$.

To finally co-ordinate such a labelling with subsequently labelling a collection of $C_4$, we need to satisfy a NTAP.
First, observe that $y$ (in $\lambda'_n(x,y)$) is not the label of any vertex if $y\neq x$. Then, in order to obtain an arithmetic progression, we might find $x+2y=3x+y$, or $y=2x$. In other words, we propose the labelling scheme $$\lambda^{(3)}_n(z)=\lambda'_n(z,2z):(z,3z,4z;5z,7z,12z;17z,19z,36z;53z,55z,\dots)$$ for collection of triangles.
For example, for a single $C_3$, we can take the labelling $(1,3,4;5,7)$, with $z=1$. 
If we consider the sequence $g(n)$ of smallest labels per cycle in this sequence of labels, assuming $z=1$, we get the recursion $g(1)=1$ and $g(k)=3g(k-1)+2$. This proves that $g(k)\in\Omega(3^k)$. Hence, the range of the labelling scheme $\lambda^{(3)}_n(z)$ grows exponentially, similar to $(\sqrt[3]{3})^n$ with the number $n$ of vertices.

We can hence use these schemes to label any collection of cycles without $C_4$, simply by interpreting the two isolates $\iota_1$ and $\iota_2$ necessary to label a collection of $\ell$ cycles (by induction hypothesis) as the first two vertices, say, $x$ and $x+y$, of the next cycle. 
Also, it is clear that 
any collection of cycles that contains at least one triangle and that is labelled this way has a NTAP.
%
Hence, we can apply \autoref{prop:induction-step+C4s} to add a collection of $C_4$ on top.
At each time, we only need two isolates for this collection of cycles, apart from one exception, when the whole collection only contains one $C_5$, where a special labelling was described in \autoref{lem:C5-C4-special}.
This proves our main theorem for 2-regular graphs.




Although our labelling strategy~$\lambda$ for 2-regular graphs~$G$ attains $\sigma(G)$, one can see that $r_\sigma(\lambda)\in\Oh(\varphi^n)$ in the worst case, where $\varphi$ is the golden ratio number. 
But if the cycle collection contains only smaller cycles, the growth rate becomes smaller. Nonetheless, it stays exponential, and it is unclear if there are labelling strategies for 2-regular graphs whose range stays polynomial.

\begin{figure}[bt]\centering
	\begin{tikzpicture}[transform shape]
 \def\xoff{4}
			\tikzset{every node/.style={circle,minimum size=0.3cm}}
			\node[draw] (a) at (0,1) {{{4}}};
			\node[draw] (b) at (-0.95106,0.30902) {{{\color{green}8}}};
			\node[draw] (c) at (-0.58779,-0.80902) {{\color{blue}{28}}};
			\node[draw] (d) at (0.58779,-0.80902) {{{36}}};
			\node[draw] (e) at (0.95106,0.30902) {{{12}}};
			\path (a) edge[-,very thick] (b);
			\path (b) edge[-,very thick] (c);
			\path (c) edge[-,very thick] (d);
			\path (d) edge[-,very thick] (e);
			\path (e) edge[-,very thick] (a);
   
   \node[draw] (a1) at (-0.9045+2*\xoff,1) {{{5}}};
			\node[draw] (b1) at (-0.9045+2*\xoff,-0.80902) {{{3}}};
			\node[draw] (c1) at (0.9045+2*\xoff,-0.80902) {{{25}}};\node[draw] (d1) at (0.9045+2*\xoff,1) {{{23}}};
			\path[draw=green] (a1) edge[-,very thick] (b1);
			\path[draw=blue] (b1) edge[-,very thick] (c1);
			\path[draw=red] (c1) edge[-,very thick] (d1);\path[draw=blue] (d1) edge[-,very thick] (a1);
  \node[draw] (a2) at (0+\xoff,1) {{{16}}};
			\node[draw] (b2) at (-0.75+\xoff,-0.80902) {{{\color{red}48}}};
			\node[draw] (c2) at (0.75+1*\xoff,-0.80902) {{{64}}};
			\path (a2) edge[-,very thick] (b2);
			\path (b2) edge[-,very thick] (c2);
			\path (c2) edge[-,very thick] (a2);
        \end{tikzpicture}
        \caption{Labelling $C_5+C_3+C_4$, with isolates $124$, $140$, and NTAP $8-28-48$.}
            \label{fig:C5C3C4}
\end{figure}
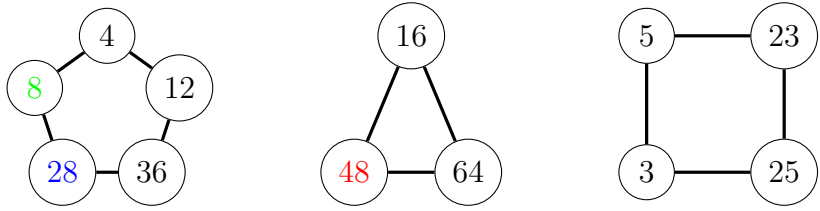

As a final comment, notice that the sequence of labellings that we propose, i.e., our labelling strategy, is not the only possible one,
as shown in \autoref{fig:C5C4C3}, where the labelling of $C_5+C_4$ as shown in \autoref{lem:C5-C4-special} is finally combined with labelling a $C_3$. Our standard labelling strategy would be a bit worse, as shown in \autoref{fig:C5C3C4}.

\begin{figure}[bt]\centering
	\begin{tikzpicture}[transform shape]
 \def\xoff{4}
			\tikzset{every node/.style={circle,minimum size=0.3cm}}
			\node[draw] (a) at (0,1) {{{2}}};
			\node[draw] (b) at (-0.95106,0.30902) {{{4}}};
			\node[draw] (c) at (-0.58779,-0.80902) {{{6}}};
			\node[draw] (d) at (0.58779,-0.80902) {{{\color{red}10}}};
			\node[draw] (e) at (0.95106,0.30902) {{{16}}};
			\path (a) edge[-,very thick] (b);
			\path (b) edge[-,very thick] (c);
			\path (c) edge[-,very thick] (d);
			\path (d) edge[-,very thick] (e);
			\path (e) edge[-,very thick] (a);
   \node[draw] (a1) at (-0.9045+\xoff,1) {{{18}}};
			\node[draw] (b1) at (-0.9045+\xoff,-0.80902) {{{26}}};
			\node[draw] (c1) at (0.9045+\xoff,-0.80902) {{{1}}};\node[draw] (d1) at (0.9045+\xoff,1) {{{9}}};
			\path[draw=green] (a1) edge[-,very thick] (b1);
			\path[draw=blue] (b1) edge[-,very thick] (c1);
			\path[draw=red] (c1) edge[-,very thick] (d1);\path[draw=blue] (d1) edge[-,very thick] (a1);
  \node[draw] (a2) at (0+2*\xoff,1) {{{\color{blue}27}}};
			\node[draw] (b2) at (-0.75+2*\xoff,-0.80902) {{{\color{green}44}}};
			\node[draw] (c2) at (0.75+2*\xoff,-0.80902) {{{71}}};
			\path (a2) edge[-,very thick] (b2);
			\path (b2) edge[-,very thick] (c2);
			\path (c2) edge[-,very thick] (a2);
        \end{tikzpicture}
        \caption{How to label $C_5+C_4+C_3$, with isolates $98$, $115$, and NTAP $10-18-26$.}
            \label{fig:C5C4C3}
\end{figure}
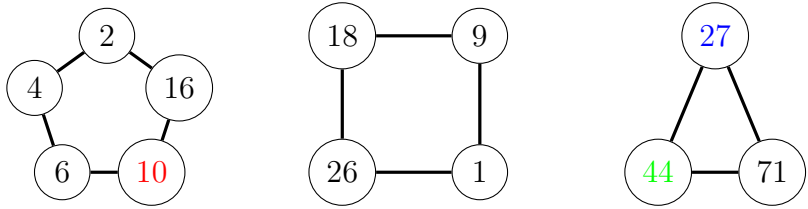

\section{Bringing paths into the game}

We are first discussing a general situation that we face after having dealt with all cycles. Here, we have to distinguish two cases: either, this cycle collection has sum number two, or it has sum number three, which means, it is a single $C_4$.

\subsection{Dealing with cycle collections of sum number two}

\begin{proposition}\label{prp:sigma2+path}
Let $G=(V,E)$ be a graph with $\sigma(G)=2$, testified by a labelling~$\lambda$ with isolate labels $\iota_1$ and $\iota_2$. If $\iota_1+\iota_2\neq \lambda(u)+\lambda(v)$ for any two vertices $u,v\in V$, then $\sigma(G+P_k)=1$ for any $k\geq 2$.  
\end{proposition}

\begin{pf}Assume $\iota_1<\iota_2$. 
Recall the Fibonacci labelling scheme for paths (\autoref{eq:fibo-scheme}). We propose to use $\lambda_{\iota_1,\iota_2}^\phi$ to label $P_k$. Clearly, $\iota_2$ and $\iota_1+\iota_2$ are the two biggest labels, labelling the second and third vertex on the path (or the isolate if $k=2$). This is an invariant that is maintained by the Fibonacci scheme: the labels of the $\ell^\text{th}$ and $(\ell+1)^\text{th}$ vertex on the path are always greater than any previous labels. By this and due to the assumption that $\iota_1+\iota_2\neq \lambda(u)+\lambda(v)$ for any two vertices $u,v\in V$, this labelling cannot introduce unwanted edges and is hence valid for $G+P$. Moreover, it will leave us with one isolate only. As $\delta(G+P)=1$, the labelling is optimal.
\qed 
\end{pf}

Notice that the argument also works if $\sigma(G)>2$; just pick the largest isolate label plus any other isolate label to produce the label $\lambda_3$ (or $\iota$ if the added path has length two).
This proves the following fact:

\begin{cor}
\label{cor:pluspath}
Let $G=(V,E)$ be a graph with $\sigma(G)\geq 2$, testified by a labelling~$\lambda$ with largest isolate labels $\iota_1$ and $\iota_2$. If $\iota_1+\iota_2\neq \lambda(u)+\lambda(v)$ for any two vertices $u,v\in V$, then $\sigma(G+P_k)\leq \sigma(G)-1$ for any $k\geq 2$. 
\end{cor}

Notice that this argument is different from the (more general) one presented in~\cite{MilRyaSmy2003} where $\sigma(G_1+G_2)\leq \sigma(G_1)+\sigma(G_2)-1$ is proved under the assumption that optimum labellings $\lambda_1$ of $G_1$ and $\lambda_2$ of $G_2$ exist such that there is an element in $\lambda_i(V(G_i))$ that is relatively prime to the largest element of $\lambda_{3-i}(V(G_{3-i}))$ for $i=1$ or $i=2$.
Also, observe that the labels will grow exponentially by a Fibonacci labelling scheme.

\subsection{Combining a 4-cycle with paths}

Given the ideas of \autoref{lem:kC4} and the results so far, one might be tempted to think that the sum number of every graph $G$ of maximum degree 2 with a disjoint copy of $C_4$ is equal to the minimum degree of $G$. Our next result shows this is not the case.

\begin{proposition}\label{prop:C4+P2}
$\sigma(C_4+P_2)=2$.
\end{proposition}

Before proving~\autoref{prop:C4+P2}, we require some observations about $C_4$. These observations also provide an indication as to why $C_4$ is different from all the other cycles (as far as the sum number is concerned, at least).

It was already shown by Harary~\cite{harary} that $\sigma(C_4)=3$. Here we present a reason for this in the following, as it indicates the way how lower bounds on $\sigma$ can be shown when the minimum degree criterion (proving $\sigma(C_r)\geq \delta(C_r)=2$ for each $r\geq 3$) is insufficient.

\begin{lem} \label{lem-4cycleexternal}
Let $C_4+G$ be a graph without isolates, and let $H$ be a sum graph of $C_4+G$. Then, all vertices corresponding to edge sums of the $C_4$ lie in $H-C_4$.
\end{lem}

In particular, this means that every sum labelling of a $C_4$ is \emph{exclusive},\footnote{This means that edge labels are among the isolate labels.} and that hence $\sigma(C_4)=3$ holds because of the next lemma (\autoref{lem-4cycle3sepsums}).

\begin{pf}[Proof of~\autoref{lem-4cycleexternal}] We will prove this by contradiction. Let the vertex labels of the $C_4$ be $(a,b,c,d)$ in cyclic order. Assume to contrary that $a+b=c$ (due to the symmetry of $C_4$, all other cases are similar). Then, we claim that all of the following are true.

\begin{enumerate}[label=(\roman*)]
    \item There is a vertex labelled $b$ in $H$.
    \item There is a vertex labelled $a+d$ in $H$.
    \item There is a vertex labelled $a+b+d$ in $H$.
\end{enumerate}

Since $b$ lies in $C_4$, (i) is true. Since $(a,d)$ is an edge and $H$ is a sum graph, (ii) is true. Finally, since $(c,d)$ is an edge, $H$ is a sum graph, and $c+d=a+b+d$, (iii) is also true. 
Now, since $H$ is a sum graph, (i), (ii), (iii) together imply that there must be an edge between the vertex labelled~$b$ and the vertex labelled $a+d$. But $b$ has exactly two neighbours, labelled~$a$ and~$c$, so $a+d$ must be one of them. We will show that either case leads to a contradiction.

If $a=a+d$, then $d=0$, which is impossible, as $H$ is a sum graph. If $c=a+d$, then $c=a+b$ implies that $b=d$, which is impossible, as $H$ is a sum graph. \qed
\end{pf}

\begin{lem} \label{lem-4cycle3sepsums}
Let $C_4+G$ be a graph without isolates, and let $H$ be a sum graph of $C_4+G$. Let $S$ be the set of numbers that correspond to the four edge sums of the $C_4$. Then, $|S|\geq 3$.
\end{lem}

\begin{pf}
We will show that $|S|\leq 2$ leads to a contradiction. That is, the four edges of the $C_4$ must have at least three~\emph{distinct} edge sums. Let the vertex labels of the $C_4$ be $(a,b,c,d)$ in cyclic order. Two edges that share a vertex cannot have the same edge sum, because then there would be two vertices with the same label. Thus, the only way that the $C_4$ can have only two distinct edge sums is if both the following hold.
\begin{align*}
a+b&=c+d\\
a+d&=c+b.
\end{align*}
Subtracting the first equation from the second, we obtain that $b-d=d-b$, or $b=d$, which is impossible in a sum graph. This completes the proof. \qed
\end{pf}

\begin{lem} \label{lem-4cycle3ariprog}
Let $C_4+G$ be a graph without isolates, and let $H$ be a sum graph of $C_4+G$. Let $S$ be the set of numbers that correspond to the four edge sums of the $C_4$. If $|S|=3$, then the three numbers in $S$ are in arithmetic progression.
\end{lem}

\begin{pf}
Let $(a,b,c,d)$ be a labelling of the $C_4$ such that $a+b=c+d$. Let
\begin{align*}
    \mathsf{sum}&=a+b=c+d;\\
    \mathsf{diff}&=c-a=b-d.
\end{align*}
We will show that the labels of the three isolates are
\begin{align*}
    \iso_1&=\mathsf{sum}-\mathsf{diff};\\
    \iso_2&=\mathsf{sum};\\
    \iso_3&=\mathsf{sum}+\mathsf{diff}.
\end{align*}
The labels of the edges $(a,b)$ and $(c,d)$ are equal to $\iso_2$, due to the definition of $\mathsf{sum}$. As for the labels of the edges $(a,d)$ and $(b,c)$, we have the following.
\begin{align*}
    a+d&=(c+d)-(c-a)=\mathsf{sum}-\mathsf{diff}=\iso_1;\\
    b+c&=(a+b)+(c-a)=\mathsf{sum}+\mathsf{diff}=\iso_3.
\end{align*}
Since $\iso_1$, $\iso_2$, $\iso_3$ are clearly in arithmetic progression, this completes the proof of~\autoref{lem-4cycle3ariprog}. \qed
\end{pf}

\noindent
Finally, we are ready to prove~\autoref{prop:C4+P2}.

\begin{pf}[Proof of~\autoref{prop:C4+P2}]
Label the $C_4$ as $(1,7,13,19)$, the $P_2$ as $(20,32)$, and the two isolates as $8$ and $44$. It is easy to check that this is a sum graph, and thus $\sigma(C_4+P_2)\leq 2$.

To prove that $\sigma(C_4+P_2)\geq 2$, assume to contrary that $\sigma(C_4+P_2)=1$. Let the labels of the vertices of the $P_2$ be $(b_1,b_2)$ (assume $b_1<b_2$), and the isolate be~$b_3$. Recall that every $C_4$ has at least three~\emph{distinct} edge labels (\autoref{lem-4cycle3sepsums}), and none of those labels can be present in the vertices of the $C_4$ itself (\autoref{lem-4cycleexternal}). Thus, the only option is that the edge labels of the $C_4$ are $b_1,b_2,b_3$. Furthermore, we know that whenever $C_4$ has exactly three edge labels, those three numbers form an arithmetic progression (\autoref{lem-4cycle3ariprog}).

Now, observe that the largest label of $P_2$ (namely, $b_2$) must be the largest label of the graph $G=C_4+P_2$, since $G$ has only one isolate. Thus, for the edge label $b_1+b_2$ of the $P_2$, we have:
\begin{equation}\label{eq:b1b2b3}
b_1+b_2=b_3.
\end{equation}
As mentioned in the previous paragraph, since $b_1 < b_2 < b_3$ are the three edge labels of the $C_4$, they are in arithmetic progression, implying that
\begin{equation}\label{eq:b3b2b1}
b_3-b_2 = b_2-b_1.
\end{equation}
With~\autoref{eq:b1b2b3} and~\autoref{eq:b3b2b1}, we get $b_2=2b_1$ and $b_3=3b_1$, or
\begin{equation}\label{eq:b_i}
b_i=ib_1\qquad\qquad\forall\ i\in\{1,2,3\}.
\end{equation}
Consider a $P_3$ subgraph of the $C_4$ such that one of the two edges of the $P_3$ is labelled $b_2$. More precisely, let the $P_3$ be $(a_1,a_2,a_3)$ such that $a_1+a_2=b_2$. Since $a_1\neq a_3$, the edge $(a_2,a_3)$ cannot be labelled~$b_2$, too. Thus, $(a_2,a_3)$ is labelled either $b_1$, or~$b_3$. That is, $a_2+a_3$ is equal to either $b_1$, or $b_3$. If $a_2+a_3=b_1$, then
\begin{align*}
a_1-a_3&=(a_1+a_2)-(a_2+a_3)\\
&=b_2-b_1\\ &=2b_1-b_1 &\text{Using~\eqref{eq:b_i}}\\ &=b_1.
\end{align*}
If $a_2+a_3=b_3$, then
\begin{align*}
a_3-a_1&=(a_2+a_3)-(a_1+a_2)\\
&=b_3-b_2\\ &=3b_1-2b_1 &\text{Using~\eqref{eq:b_i}}\\ &=b_1.
\end{align*}
Therefore, either $a_1=b_1+a_3$, or $a_3=b_1+a_1$. In other words, either $(b_1,a_3)$ is an edge, or $(b_1,a_1)$ is an edge. In either case, there is an edge between the $C_4$ and the $P_2$, which is a contradiction because they are supposed to be disjoint. \qed
\end{pf}

\noindent
On the other side, we can prove:
\begin{lem}
$\sigma(C_4+P_k)=1$ for all $k\geq 3$.
\end{lem}

\begin{pf}
If $k=3$, label the $C_4$ as $\{1,3,9,11\}$ and the $P_3$ as $\{12,4,16\}$, with the isolate being~$20$. If $k\geq 4$, then the first four labels of the path $P_k=(v_1,v_2\dots,v_k)$ are $\lambda(v_1)=12$, $\lambda(v_2)=4$, $\lambda(v_3)=16$, $\lambda(v_4)=20$. After that, we simply continue in the Fibonacci fashion, i.e., $\lambda(v_{i+1})=\lambda(v_i)+\lambda(v_{i-1})$, with the label of the isolate being $\lambda(v_k)+\lambda(v_{k-1})$. It is easy to check that no unwanted edges are introduced.
\qed
\end{pf}

The general algebraic strategy can be best seen by labelling the $C_4$ with $\{1,3,9,11\}$, with $a<b$ being the smallest numbers. We assume that $a+d=b+c$, i.e., $d=b+c-a$. Moreover, we label the three path vertices with $\{a+d,a+b,(a+d)+(a+b)=a+2b+c\}$.
In order to save on isolates, we also require that $c+d=b+2c-a$ equals $(a+2b+c)+(a+b)=2a+3b+c$, which implies $c=3a+2b$.
In summary, given small numbers $a<b$, we construct the further labels of $C_4$ as $3a+2b$ and $2a+3b$. Then, the labels on the path would be $3(a+b)$, $a+b$, $4(a+b)$, with the isolate $5(a+b)$.

If we want to label $C_4+P_k$ with $k\geq 4$, it is possible to save on the size of the labels by starting with labelling the $C_4$ with $\{1,2,6,11\}$ and the $P_4$ with $\{17,3,8,12\}$, plus one isolate labeled~$20$.
Further savings are possible if we label the $C_4$ with $\{2,5,8,11\}$, as done as a standard throughout this paper.
The first five labels of the path $P_k$, $k\geq 5$, with vertices $v_1,\dots,v_k$ are then:
$\lambda(v_1)=26$, $\lambda(v_2)=13$, $\lambda(v_3)=7$, $\lambda(v_4)=19$,
$\lambda(v_5)=20$. If $k=5$, then $39$ would be the label of the isolate. Otherwise, we just continue in a Fibonacci-style, i.e., $\lambda(v_{i+1})=\lambda(v_i)+\lambda(v_{i-1})$, with $\lambda(v_k)+\lambda(v_{k-1})$ being the isolate. Again, no unwanted edges are introduced.

There is only one case left over to complete the picture:
\begin{lem}\label{lem:C4+2P2}
$\sigma(C_4+2P_2)=1$.
\end{lem}

\begin{pf}
By \autoref{prop:C4+P2}, $\sigma(C_4+P_2)=2$. The labelling satisfies the requirements of \autoref{prp:sigma2+path}, which shows the claim.\qed 
\end{pf}

\subsection{Combining cycles with more than one path}

\noindent
The following proposition also covers the case of pure path collections.

\begin{proposition}\label{prp:sigma1+paths}
Let $G$ be a graph with $\sigma(G)=1$.
Then, $\sigma(G+P_k)=1$ for any  $k\geq 2$.
\end{proposition}

\begin{pf}
Let $\iota$ be the label of the isolate of a labelling $\lambda$ of $G$ certifying its sum number to be~1. Then, $\iota$ is bigger than any vertex label of~$G$. Therefore, labelling $P_k$ with the Fibonacci labelling scheme $\lambda_{\iota,2\iota}^\phi$, as introduced in \autoref{eq:fibo-scheme} in general form, labels $G+P_k$ (together with $\lambda$) with only one isolate, not creating conflicts, as all edge labels of $G$ are smaller than $2\iota$.\qed 
\end{pf}

We already saw (or will see soon) that a cycle collection plus one path has sum number one with one exception, which is $C_4+P_2$.
As we will fix the only remaining case of $C_4+2P_2$ separately in \autoref{lem:C4+2P2}, we can conclude:

\begin{proposition}
Let $C$ be a collection of cycles and $P$ be a collection of at least two paths. Then, $\sigma(C+P)=1$.
\end{proposition}

\begin{pf}
Except for the case of $C_4+P_2$, we know that $\sigma(C+P_k)=1$ for any collection of cycles $C$ and any $k\geq 3$. As we consider paths in decreasing length, we will pick a cycle of length at least three to be considered first if there is any. Therefore, if $P$ is a collection of at least two paths, then we can conclude  $\sigma(C+P)=1$ either by \autoref{prp:sigma1+paths} directly, or by first taking \autoref{lem:C4+2P2}.\qed 
\end{pf}

\section{Conclusion and Open Problems}

We have explained that the labelling of a $C_5$ as proposed in \cite{SinTiwTri2021} is not working correctly. This leaves the spum-minimization problem open for this particular small graph.
But this question easily generalizes to nearly all graphs with maximum degree of at most two as discussed in this paper. In most cases, we only found labellings with labels of exponential size. This might be 
necessary, but for such statements, we do not have any proof idea.

Our main result concerns the sum number of (all) graphs of maximum degree two. Kratochv\'il, Miller and Nguyen posed in~\cite{KraMilNgu2001} two conjectures that are tightly related to our paper; we will formulate them as questions below.
    \begin{itemize}
        \item Given two graphs $G_1, G_2$ with $\sigma(G_1)=\sigma(G_2)=1$, is it true that the sum number of their disjoint union is always one?
        \item More generally: given two graphs $G_1, G_2$, is it true that $\sigma(G_1+G_2)\leq \sigma(G_1)+\sigma(G_2)-1$?
    \end{itemize}
    Observe that we did resolve the first question if $G_2$ is a path (\autoref{prp:sigma1+paths}), but the general question is still open. Upon some thought, it can be seen that the general question is related to the following natural combinatorial question: Find a characterization of all graphs with sum number one, also known as~\emph{unit graphs} in the literature. We also refer to a recent paper~\cite{KonKNPST2018} that studies variations of this question. Finally, a slightly weaker but more structured notion of sum labelling (called arithmetic graphs) could lend some ideas that might help in resolving this question (and more optimistically, the second question)~\cite{AchHeg89}.
    
    Apart from these combinatorial questions, the basic complexity questions concerning the graph parameter $\sigma$ and $r_\sigma$ are open. For instance, is it \textsf{NP}-hard to decide if, given a graph~$G$ and a number $k$, $\sigma(G)\leq k$ holds? One of our own motivations to study graphs of maximum degree two was to see if one could use the operation of graph union to piece gadgets together for this and similar questions. But we are still far from this, as even these seemingly easy questions concerning graphs of maximum degree two are non-trivial to solve.

\bibliographystyle{elsarticle-num}
\bibliography{biblio.bib}

\end{document}